\documentclass[letterpaper,10pt,conference]{ieeeconf}
\IEEEoverridecommandlockouts
\DeclareUnicodeCharacter{FF0E}{{$\mathcal N$}}
\usepackage{graphicx}
\usepackage{macros_conf}
\usepackage{amsfonts}
\usepackage{nccmath}
\usepackage{amsmath}
\usepackage{amssymb}
\usepackage{arydshln}
\usepackage{fleqn}
\usepackage{here}
\usepackage{macros_e}
\usepackage{mathfonts}
\usepackage{footnote}

\everymath={\displaystyle}
\title{
\LARGE \bf 
Stability Analysis of Recurrent Neural Networks by\\ IQC with Copositive Mutipliers
}
\author{Yoshio Ebihara, Hayato Waki, Victor Magron, \\
Ngoc Hoang Anh Mai, Dimitri Peaucelle, Sophie Tarbouriech
\thanks{
Y. Ebihara is with the 
Graduate School of Information Science and 
Electrical Engineering, Kyushu University, 
744 Motooka, Nishi-ku, Fukuoka 819-0395, Japan, 
he was also with 
LAAS-CNRS, Universit\'{e} de Toulouse, CNRS, Toulouse, France,   
in 2011.  
H. Waki is with the Institute of Mathematics for Industry,
Kyushu University, 744 Motooka, Nishi-ku, Fukuoka 819-0395, Japan.  
V. Magron, N. H. Mai, D. Peaucelle, and S. Tarbouriech are with 
LAAS-CNRS, Universit\'{e} de Toulouse, CNRS, Toulouse, France.  
}%
}

\begin{document}
\maketitle
\thispagestyle{empty}
\pagestyle{empty}

\begin{abstract}
 This paper is concerned with the stability analysis of the
 recurrent neural networks (RNNs) by means of the
 integral quadratic constraint (IQC) framework.
 The rectified linear unit (ReLU) is typically employed 
 as the activation function of the RNN,
 and the ReLU has specific nonnegativity properties
 regarding its input and output signals.
 Therefore, it is effective if we can derive IQC-based stability
 conditions
 with multipliers taking care of such nonnegativity properties.  
 However, such nonnegativity (linear) properties are hardly captured
 by the existing multipliers defined on the positive semidefinite cone.
 To get around this difficulty,
 we loosen the standard positive semidefinite cone 
to the copositive  cone,  
 and employ copositive multipliers to capture the nonnegativity
 properties.
 We show that, within the framework of the IQC, 
 we can employ copositive multipliers (or their inner approximation)
 together with existing multipliers such as
 Zames-Falb multipliers and polytopic bounding multipliers, 
 and this directly enables us to ensure that the introduction of 
 the copositive multipliers leads to better (no more conservative) results.  
 We finally illustrate the effectiveness of 
 the IQC-based stability conditions with the copositive multipliers 
 by numerical examples.  

\noindent
 {\bf Keywords: recurrent neural networks,
 rectified linear units, stability, IQC, nonnegative signals, 
copositive multipliers.  
}
\end{abstract}


\section{Introduction}
\label{sec:intro}

A recurrent neural network (RNN) is a class of deep neural networks 
and able to imitate the behavior of dynamical systems
due to its feedback mechanism.  
The effectiveness of the RNN is widely recognized in speech recognition, 
natural language processing, and image recognition 
\cite{Barabanov_IEEENN2002,Zhang_IEEENN2014,Salehinejad_2018}.  
Even though new architectures such as transformer \cite{Vaswani_NIPS2017}
have been developed recently, it is expected 
that the RNN retains its position 
as one of the fundamental and important elements in deep neural networks.  

Even though the feedback mechanism is the key of the RNN
and distinguishes the RNN from other feedforward networks,   
the existence of the feedback mechanism may cause network instability.  
Therefore the stability analysis of the RNN has been an important issue
in the machine learning community 
\cite{Barabanov_IEEENN2002,Zhang_IEEENN2014,Salehinejad_2018}.  
From control theoretic viewpoint, 
we can readily apply the small gain theorem \cite{Khalil_2002}
to the stability analysis of a given RNN by representing it 
as a feedback connection with a linear time-invariant (LTI) system  and
a static nonlinear activation function  typically being 
a rectified linear unit (ReLU) for the RNN.  
It is nonetheless true that the standard small gain theorem leads to 
conservative results 
since it does not take into account the important property that 
the ReLU returns only nonnegative signals.  
This motivated us to analyze the $l_2$ induced norm of LTI systems
for nonnegative input signals in \cite{Ebihara_EJC2021}, 
which is referred to as the $l_{2+}$ induced norm in this paper.  
We characterized an upper bound of 
the $l_{2+}$ induced norm 
by copositive programming \cite{Dur_2010}, and then
derived a numerically tractable semidefinite program (SDP)
for (in general loosened) upper bound computation. 
We finally derived an  $l_{2+}$-induced-norm-based 
(scaled) small gain theorem for the stability analysis of the RNN and 
illustrated its effectiveness by numerical examples.  

We believe that the treatments in \cite{Ebihara_EJC2021}
brought some new insights for the stability analysis of
feedback systems constructed from LTI systems and 
nonlinear elements (i.e., Lur'e systems).  
However, the $l_{2+}$-induced-norm-based (scaled) small gain condition
might be shallow in view of 
the advanced integral quadratic constraint (IQC) theory \cite{Megretski_IEEE1997}.  
We acknowledge the fact that, for the stability analysis of Lur'e systems, 
the effectiveness of 
the IQC-based approaches with Zames-Falb multipliers \cite{Zames_SIAM1968}
are widely recognized, see, e.g., \cite{Fetzler_IFAC2017,Fetzler_IJRN2017}.  
Therefore it is strongly preferable if we can build the nonnegativity-based approach
upon the powerful IQC-based framework.  
Such an extension seems hard, since, as the denomination IQC says, 
the existing multipliers capture the properties of nonlinear elements 
with {\it quadratic} constraints on their input-output signals, 
whereas the nonnegativity property of the RNN (i.e., ReLU) is essentially 
{\it linear} constraints on the input-output signals.   
To get around this difficulty, 
we loosen the standard positive semidefinite cone to the copositive cone 
and employ copositive multipliers to handle 
the linear (nonnegativity) constraints on the input-output signals of the RNN.   
As clarified later on, this can be done in such a sound way that 
the proposed IQC-based stability condition with the copositive multipliers
encompasses the results in \cite{Ebihara_EJC2021} as particular cases.  
Then, by applying an inner approximation to the copositive cone, 
we derive numerically tractable IQC-based SDPs for the stability analysis of the RNN.  
We show that, within the framework of IQC, 
we can employ copositive multipliers (or their inner approximation)
together with existing multipliers such as
the Zames-Falb multipliers and polytopic bounding multipliers, 
and this directly enables us to ensure that the introduction of 
the copositive multipliers leads to better (no more conservative) results.  
We finally illustrate the effectiveness of 
the IQC-based stability conditions with the copositive multipliers 
by using the same numerical examples as in \cite{Ebihara_EJC2021}.
Related works include \cite{Anderson_IEEENN2007,Yin_IEEE2021,Fazlyab_IEEE2021,Revay_LCSS2021}, but
again the novel contribution of the present paper is
capturing the behavior of ReLUs by 
copositive multipliers within the framework of IQCs.    

Notation: 
The set of $n\times m$ real matrices is denoted by $\bbR^{n\times m}$, and
the set of $n\times m$ entrywise nonnegative matrices is denoted
by $\bbR_+^{n\times m}$. 
For a matrix $A$, we also write $A\geq 0$ to denote that  
$A$ is entrywise nonnegative. 
We denote the set of $n\times n$ real symmetric matrices by $\bbS^n$.
For $A\in\bbS^n$, we write $A\succ 0\ (A\prec 0)$ to
denote that $A$ is positive (negative) definite.  
For $A\in\bbR^{n\times n}$, we define $\He\{A\}:=A+A^T$.
For $A\in\bbR^{n\times n}$ and $B\in\bbR^{n\times m}$,
$(\ast)^TAB$ is a shorthand notation of $B^TAB$.  
We denote by $\bbD_{++}^n\subset\bbR^{n\times n}$ the set of diagonal matrices
with strictly positive diagonal entries.  
In addition, we denote by $\bbD^n[\alpha,\beta]$
the set of diagonal matrices whose diagonal entries are all within 
the closed interval $[\alpha,\beta]$.  
Moreover, $\bbD_\mathrm{ver}^n[\alpha,\beta]\subset \bbD^n[\alpha,\beta]$
is the set of $2^n$ matrices corresponding to the vertices of $\bbD^n[\alpha,\beta]$.  
A matrix $M\in\bbR^{n\times n}$ is said to be 
Z-matrix if $M_{ij}\le 0$ for all $i\neq j$.  
Moreover, $M$ is said to be doubly hyperdominant if
it is a Z-matrix and 
$M \one_n \ge 0$, $\one_n^T M \ge 0$,   
where $\one_n\in\bbR^n$ stands for the all-ones-vector.  
In this paper we denote by $\DHD^n\subset \bbR^{n\times n}$ the set of 
doubly hyperdominant matrices.  

For the discrete-time signal $w$ defined 
over the time interval $[0,\infty)$, we define
\[
\begin{array}{@{}l}
 \|w\|_{2}:=\sqrt{\sum_{k=0}^{\infty}|w(k)|_2^2}
\end{array}
\]
where for $v\in\bbR^{n_v}$ we define
$|v|_2:=\sqrt{\sum_{j=1}^{n_v} v_j^2}$.  
We also define
\[
\begin{array}{@{}l}
 l_{2}  :=\left\{w:\ \|w\|_{2}<\infty \right\},\\ 
 l_{2+} :=\left\{w:\ w\in l_{2},\ w(k)\ge 0\ (\forall k\ge 0)\right\}  
\end{array}
\]
and
\[
\begin{array}{@{}l}
 l_{2e}  :=\left\{w:\ w_\tau\in l_2,\ \forall \tau\in[0,\infty) \right\}
\end{array}
\]
where $w_\tau$ is the truncation of the signal $w$ up to the time instant
$\tau$
and defined by
\[
 w_\tau(k)=
 \left\{
 \begin{array}{cc}
  w(k)& (k\le \tau), \\
  0& (k> \tau).  \\
 \end{array}
 \right.  
\]
For an operator $H:\ l_{2e}\ni w \to z \in l_{2e}$, 
we define its (standard) $l_2$ induced norm by
\begin{equation}
\|H\|_{2}:=\sup_{w\in l_2,\ \|w\|_2=1} \ \|z\|_2.  
\label{eq:l2norm}
\end{equation}
We also define
\begin{equation}
\|H\|_{2+}:=\sup_{w\in l_{2+},\ \|w\|_2=1} \ \|z\|_2.  
\label{eq:l2+norm}
\end{equation}
This is a variant of the $l_2$ induced norm introduced in \cite{Ebihara_EJC2021}
and referred to as the $l_{2+}$ induced norm in this paper.  
We can readily see that $\|H\|_{2+}\le \|H\|_{2}$.  

\section{Copositive Programming}
\label{sec:cop}

Copositive programming (COP) is a convex optimization problem in which 
we minimize a linear
objective function over the 
linear matrix inequality (LMI) constraints on the copositive cone   
\cite{Dur_2010}.  
In this section, we summarize its basics.  

\subsection{Convex Cones Related to COP}

Let us review the definition and the property of convex
cones related to COP.
\begin{definition}\cite{Berman_2003}\label{conedef} 
The definition of proper cones 
$\PSD_n$, $\COP_n$, $\CP_n$, $\NN_n$, and $\DNN_n$ in 
$\bbS^n$ are as follows.
\begin{enumerate}
  \item 
  $\PSD_n:=\{P\in\bbS^n:\ \forall x\in\bbR^n,\
	x^TPx\geq 0\}=\{P\in\bbS^n:\ \exists B\
	\mbox{s.t.}\ P=BB^T\}$ is called {\it the positive semidefinite cone}.
  \item
  $\COP_n:=\{P\in\bbS^n:\ \forall x\in\bbR_{+}^n,\
	x^TPx\geq 0\}$ is called {\it the copositive cone}.
  \item
  $\CP_n:=\{P\in\bbS^n:\ \exists B\ge 0\
	\mbox{s.t.}\ P=BB^T\}$ is called {\it the completely positive cone}.
  \item
  $\NN_n:=
  \{P\in\bbS^n:\ P\geq 0\}
  $
  is called {\it the nonnegative cone}.
  \item
  $\PSD_n+\NN_n:=\{P+Q:\ P\in\PSD_n,\
	Q\in\NN_n\}$．
       This is the Minkowski sum of the positive semidefinite cone and 
       the nonnegative cone.
  \item
  $\DNN_n:=\PSD_n\cap\NN_n$ is called 
  {\it the doubly nonnegative cone}.
\end{enumerate}
\end{definition}

From Definition \ref{conedef}, 
we clearly see that the following inclusion relationships hold: 
\begin{equation}
 \CP_n\subset\DNN_n\subset\PSD_n\subset\PSD_n+\NN_n\subset\COP_n,
\label{eq:inc1}
\end{equation}
\begin{equation}
 \CP_n\subset\DNN_n\subset\NN_n\subset\PSD_n+\NN_n\subset\COP_n.
\label{eq:inc2}
\end{equation}
In particular, when $n\leq 4$, it is known that 
$\COP_n=\PSD_n
+\NN_n$ and $\CP_n =\DNN_n$ hold \cite{Berman_2003}．
On the other hand, as for the duality of these cones, 
$\COP_n$ and $\CP_n$ are dual to each other, 
$\PSD_n+\NN_n$ and $\DNN_n$ are dual to each other, 
and $\PSD_n$ and $\NN_n$ are self-dual.
It is also well known that the interior of the cone $\PSD_n$ can be
characterized by
\[
\begin{array}{@{}lcl}
  \PSD_n^\circ&=&\{P\in\bbS^n:\  \forall
   x\in\bbR^n\backslash\{0\},\ x^TPx>0\}\\
 &=&
  \{P\in\bbS^n:\  \exists B\ \mbox{s.t.}\ P=BB^T,\
       \mathrm{rank}(B)=n\}.  
\end{array}
\]
%

\subsection{Basic Properties of COP}

COP is a convex optimization problem on the copositive cone and
its dual is a convex optimization problem on the completely positive cone.  
As mentioned in \cite{Dur_2010}, 
the problem to determine whether a given symmetric matrix is
copositive or not is a co-NP complete problem, and
the problem to determine whether a given symmetric matrix is
completely positive or not is an NP-hard problem.  
Therefore, it is hard to solve COP 
numerically in general.  
However, since the problem to determine whether a given matrix is in $\PSD+\NN$
or in $\DNN$ can readily be reduced to SDPs, 
we can numerically solve the convex optimization problems on the cones
$\PSD+\NN$ and $\DNN$ easily.  
Moreover, when $n\leq 4$, it is known that $\COP_n=\PSD_n+\NN_n$ and 
$\CP_n=\DNN_n$ as stated above, 
and hence those COPs with $n\leq 4$ can be reduced to SDPs.  

\section{IQC-Based Stability Analysis of RNN with ReLU}
\label{sec:RNN}

\subsection{Basics of RNN and Stability}
\label{sub:RNN}

Let us consider the dynamics of 
the discrete-time RNNs typically described by
\begin{equation}
G:\ \left\{
\arraycolsep=0.5mm
 \begin{array}{ccl}
  x(k+1)&=&\Lambda x(k)+\Win w(k)+v(k),\\
  z(k)&=&\Wout x(k),\\
  w(k)&=&\Phi(z(k)+s(k)).  
 \end{array}
\right.  
\label{eq:RNN}
\end{equation}
Here $x\in\bbR^n$ is the state and 
$\Lambda\in\bbR^{n\times n}$, 
$\Wout\in\bbR^{m\times n}$, 
$\Win\in\bbR^{n\times m}$ are constant matrices
with $\Lambda$ being Schur-Cohn stable. 
We assume $x(0)=0$.  
On the other hand, note that 
$s:\ [0,\infty)\to \bbR^m$ and $v:\ [0,\infty)\to \bbR^n$ are 
external input signals  
and $\Phi:\ \bbR^m\to \bbR^m$ is the
static activation function typically being nonlinear.  
The matrices $\Win$ and $\Wout$ are constructed from 
the weightings of the edges in RNN.  

In this paper, we consider the typical case where the activation function is
the (entrywise)
rectified linear unit (ReLU) whose input-output property is given by
\begin{equation}
\begin{array}{@{}l}
 \Phi(\xi)=\left[\ \phi(\xi_1)\ \cdots\ \phi(\xi_m)\ \right]^T,\\
 \phi: \bbR\to \bbR,\quad
 \phi(\eta)=
  \left\{
   \begin{array}{cc}
    \eta & (\eta\ge 0),  \\
    0 & (\eta< 0).  \\
   \end{array}
  \right.    
\end{array}  
\label{eq:ReLU}
\end{equation}
We can readily see that $\|\Phi\|_2=1$.  
It should be noted that 
the system $G_0$
essentially makes the feedback loop with
the ReLU $\Phi$ where
\begin{equation}
 G_0:=
 \left[
 \begin{array}{c|c}
  \Lambda & \Win \\ \hline
  \Wout & 0 
 \end{array}
 \right].  
\label{eq:G0}
\end{equation}
%

Since here we are dealing with nonlinear systems, 
it is of prime importance to clarify the definition of 
``stability.''
The definition we employ for the analysis of RNN is as follows.  
\begin{definition}\cite{Khalil_2002} (Finite Gain $l_2$ Stability) 
An operator $H:\  l_{2e}\ni u \to y\in l_{2e}$
is said to be finite gain $l_2$ stable
if there exists a nonnegative constant $\gamma$ such that
$\|y_\tau\|_2 \le \gamma \|u_\tau\|_2$ holds
for any $u \in l_{2e}$ and $\tau \in [0,\infty)$. 
\end{definition}

In the following, we analyze the 
finite gain $l_2$ stability of the operator in RNN
with respect to
the input $[\ s^T\ v^T\ ]^T\in l_{2e}$
and the output $[\ z^T\ w^T\ ]^T \in l_{2e}$.  
Note that the feedback connection in the RNN 
is well-posed since its dynamics is given by the state-space equation
\rec{eq:RNN}.  
We also note that we implicitly use the causality of $G$ and $\Phi$
in the following.  

\subsection{IQC-Based Basic Stability Condition}

It is known that the framework of 
Integral Quadratic Constraint (IQC) \cite{Megretski_IEEE1997} is helpful
in capturing the nonlinearity in feedback systems and
obtaining less conservative results for stability analysis.  
The basic IQC-based stability condition 
for RNN with ReLU can be summarized by the next theorem.  
\begin{theorem}
For any input signal $\xi\in\l_{2e}$ and 
output signal $\zeta\in l_{2e}$
of ReLU $\Phi$ such that $\zeta=\Phi \xi$, 
suppose $\Pi\in\bbS^{2m}$ satisfies 
the time-domain 
(discrete-time version of) IQC given by
\begin{equation}
 \sum_{k=0}^\tau
 \left[
 \begin{array}{c}
 \xi(k) \\
 \zeta(k) \\
 \end{array}
 \right]^T
 \Pi
 \left[
 \begin{array}{c}
 \xi(k) \\
 \zeta(k) \\
 \end{array}
 \right]\ge 0
\label{eq:IQC0}
\end{equation}
for any $\tau\in[0,\infty)$.  
Then, the RNN given by \rec{eq:RNN}
with ReLU $\Phi$ given by \rec{eq:ReLU}
is finite-gain $l_2$ stable if there exist
$P\in\PSD_n$ and 
$S\in\bbD_{++}^m$ such that
\begin{equation}
\scalebox{0.95}{$
\begin{array}{@{}l}
\left[
\begin{array}{cc}
 -P & 0 \\
 0 & -S
\end{array}
\right]+
 \left[
 \begin{array}{cc}
  \Lambda & \Win \\
  \Wout & 0\\
 \end{array}
 \right]^T
 \left[
 \begin{array}{cc}
  P & 0\\
  0 & S\\
 \end{array}
 \right]
 \left[
 \begin{array}{cc}
  \Lambda & \Win \\
  \Wout & 0\\
 \end{array}
 \right]\\
\hspace*{5mm}
+
\left[
 \begin{array}{cc}
  \Wout & 0 \\
  0 & I_m \\
 \end{array}
\right]^T\Pi
\left[
 \begin{array}{cc}
  \Wout & 0 \\
  0 & I_m \\
 \end{array}
\right]
\prec 0.  
\end{array}$}
\label{eq:RNNIQC} 
\end{equation}
\label{th:IQC}
\end{theorem}
\begin{proofof}{\rth{th:IQC}}
Suppose \rec{eq:RNNIQC} holds with 
$P=\hatP\in \PSD_n$ and  
$S=\hatS\in \bbD_{++}^m$.  
Then, it is very clear that there exist 
$\varepsilon>0$ and $\nu>0$ such that
\[
\scalebox{0.95}{$
\begin{array}{@{}l}
 M(\hatP,\varepsilon,\nu)\prec 0,\ \mbox{with}\\  
M(\hatP,\varepsilon,\nu):=
\arraycolsep=0.3mm
\left[
\begin{array}{cccc}
 -\hatP+\varepsilon^2 \Wout^T\Wout & 0 & 0 & 0 \\
 0 & -\hatS & 0 & 0 \\ 
 0 & 0 & -\nu^2 I_n& 0 \\ 
 0 & 0 & 0 &-\nu^2 I_m \\ 
\end{array}
\right]\\
\arraycolsep=0.3mm
+(\ast)^T
 \left[
 \begin{array}{cc}
  \hatP & 0\\
  0 & \hatS\\
 \end{array}
 \right]
 \left[
 \begin{array}{cccc}
  \Lambda & \Win & I_n & 0 \\
  \Wout & 0  & 0  & I_m\\
 \end{array}
 \right]
+(\ast)^T
  \Pi
 \left[
 \begin{array}{cccc}
  \Wout & 0   & 0  & I_m \\
  0     & I_m & 0  & 0 \\
 \end{array}
 \right].  
\end{array}$}
\]
Then, along the trajectory of the RNN 
for the input signals $v\in l_{2e}$ and $s\in l_{2e}$, 
we have
\[
 \left[
 \begin{array}{c}
  x(k)\\
  w(k)\\
  v(k)\\
  s(k)\\
 \end{array}
 \right]^T
M(\hatP,\varepsilon,\nu) 
 \left[
 \begin{array}{c}
  x(k)\\
  w(k)\\
  v(k)\\
  s(k)\\
 \end{array}
 \right]
 \le 0\ (k=0,1,\cdots)  
\]
or equivalently, 
\[
\begin{array}{@{}l}
\varepsilon^2 z(k)^Tz(k)-x(k)^T\hatP x(k)+x(k+1)^T\hatP x(k+1)\\
+(z(k)+s(k))^T\hatS(z(k)+s(k))-w(k)^T\hatS w(k)\\
-\nu^2 v(k)^Tv(k)-\nu^2 s(k)^Ts(k)\\
+\left[
 \begin{array}{c}
 z(k)+s(k) \\
 w(k) \\
 \end{array}
 \right]^T
 \Pi
 \left[
 \begin{array}{c}
 z(k)+s(k) \\
 w(k) \\
 \end{array}
 \right] \le 0\\
(k=0,1,\cdots).  
\end{array}
\]
Here, since $\|\Phi\|_2 = 1$ and $\hatS\in\bbD_{++}^m$, we have
\[
(z(k)+s(k))^T\hatS (z(k)+s(k))-w(k)^T\hatS w(k)\ge 0 
\]
and hence 
\[
\begin{array}{@{}l}
\varepsilon^2 z(k)^Tz(k)-x(k)^T\hatP x(k)+x(k+1)^T\hatP x(k+1)\\
-\nu^2 v(k)^Tv(k)-\nu^2 s(k)^Ts(k)\\
+
 \left[
 \begin{array}{c}
 z(k)+s(k) \\
 w(k) \\
 \end{array}
 \right]^T
 \Pi
 \left[
 \begin{array}{c}
 z(k)+s(k) \\
 w(k) \\
 \end{array}
 \right] \le 0\\
(k=0,1,\cdots).    
\end{array}
\]
By summing up the above inequality up to $k=\tau$, 
we have
 \[
 \scalebox{0.82}{$
 \begin{array}{@{}l}
 x(\tau+1)^T\hatP x(\tau+1)
+\varepsilon^2\sum_{k=0}^\tau |z(k)|_2^2
-\nu^2 \left(\sum_{k=0}^\tau |v(k)|_2^2+ \sum_{k=0}^\tau |s(k)|_2^2\right)\\
  \renewcommand{\arraystretch}{0.90}
+ \sum_{k=0}^\tau
 \left[
 \begin{array}{c}
 z(k)+s(k) \\
 w(k) \\
 \end{array}
 \right]^T
 \Pi
 \left[
 \begin{array}{c}
 z(k)+s(k) \\
 w(k) \\
 \end{array}
 \right] \le 0.  
 \end{array}$}
\]
Since $\hatP\in\PSD_n$ and since \rec{eq:IQC0} holds, 
we can readily conclude from the above inequality that
\[
 \|z_\tau\|_2^2 \le \frac{\nu^2}{\varepsilon^2}
 \left(\|v_\tau\|_2^2+ \|s_\tau\|_2^2  \right)
\]
or equivalently, 
\[
  \renewcommand{\arraystretch}{0.90}
 \|z_\tau\|_2 \le \frac{\nu}{\varepsilon} 
\left\|
\left[
\begin{array}{c}
 s_\tau\\
 v_\tau
\end{array}\right]
\right\|_2.  
\]
With this inequality and 
\[
 \renewcommand{\arraystretch}{0.95}
\begin{array}{@{}lcl}
 \|w_\tau\|_{2} &=& \|(\Phi (z+s))_\tau\|_2\\
 &=&  \|\Phi (z+s)_\tau\|_2\\
 &\le&  \|\Phi\|_2 (\|z_\tau\|_2+\|s_\tau\|_2)\\
 &=&  \|z_\tau\|_2+\|s_\tau\|_2,   
\end{array}
\]
we arrive at the conclusion that 
\[
\begin{array}{@{}lcl}
  \renewcommand{\arraystretch}{0.90}
\left\|
\left[
\begin{array}{c}
 z_\tau\\
 w_\tau
\end{array}\right]
\right\|_2\le 
\sqrt{\frac{\nu^2}{\varepsilon^2}+2} \left\|
\left[
\begin{array}{c}
 s_\tau\\
 v_\tau
\end{array}\right]
\right\|_2  
\end{array}
\]
holds for any $v\in l_{2e}$, $s\in l_{2e}$ and $\tau\in[0,\infty)$.  
This completes the proof.  
\end{proofof}
\begin{remark}
Since $G_0$ defined in \rec{eq:G0} 
makes the feedback loop with $\Phi$, 
and since $\| \Phi\|_2=1$, 
it is very clear that 
the small gain condition $\|G_0\|_2<1$ is a 
sufficient condition for the stability of the RNN with the ReLU.  
In addition, it is not hard to see that the ReLU $\Phi$ satisfies
$\Phi(\xi)=(D^{-1}\Phi D)(\xi)$ for any 
$D\in\bbD_{++}^{n}$.  
Therefore the scaled small gain condition
$\|D^{-1}G_0D\|_{2}< 1$ with $D\in\bbD_{++}^{n}$
is also a sufficient condition for the stability.  
It should be noted that \rec{eq:RNNIQC} with $\Pi=0$
corresponds to the scaled small gain condition, 
and that \rec{eq:RNNIQC} with $\Pi=0$ and $S=I_m$
corresponds to the small gain condition \cite{Khalil_2002}.
In this sense, the IQC-based stability condition in
\rth{th:IQC} encompasses these basic stability conditions.  
\end{remark}
%

\section{Concrete Multipliers Capturing the Properties of ReLU}

\subsection{Zames-Falb Multiplier}

In this section, we summarize the arguments of
\cite{Fetzler_IFAC2017} on the discrete-time 
Zames-Falb multipliers \cite{Zames_SIAM1968}.  
By following \cite{Fetzler_IFAC2017}, we first introduce the following definitions.  
\begin{definition}\cite{Fetzler_IFAC2017} 
Let $\mu\le 0 \le \nu$.  
Then the nonlinearity $\phi:\ \bbR\to\bbR$ is slope-restricted, 
in short $\phi \in \slope(\mu,\nu)$, if $\phi(0)=0$ and
\[
  \mu\leq \dfrac{\phi(x)-\phi(y)}{x-y}
 \le
  \sup_{x\neq y}\dfrac{\phi(x)-\phi(y)}{x-y}
  <\nu
\]
 for all $x,y\in\bbR$, $x\neq y$.
 On the other hand, the nonlinearity $\phi$ is said to be sector-bounded if 
\[
  (\phi(x)-\alpha x)(\phi(x)-\beta x)\leq 0
  \ (\forall x\in\bbR)
\]
for some $\alpha \le 0 \le \beta$.  
This is expressed as $\phi\in\sec[\alpha,\beta]$.  
\end{definition}

The main result of \cite{Fetzler_IFAC2017}
on the discrete-time Zames-Falb multipliers 
for slope-restricted nonlinearities can be
summarized by the next lemma.  
\begin{lemma}\cite{Fetzler_IFAC2017}
For a given nonlinearity $\phi\in \slope(\mu,\nu)$ with 
$\mu\le 0\le \nu$, 
let us define $\Phi:\ \bbR^m\to \bbR^m$
by the first equation in \rec{eq:ReLU}.  
Assume $M\in\DHD^m$.  Then we have
\[
\scalebox{0.93}{$
\begin{array}{@{}l}
 (\ast)^T
\left[
 \begin{array}{cc}
  0 & M^T\\
  M & 0 
 \end{array}
\right]
\left(
 \left[
 \begin{array}{cc}
  \nu I_m & -I_m\\
  -\mu I_m & I_m\\
 \end{array}
 \right]
 \left[
 \begin{array}{c}
  x\\
  \Phi(x)
 \end{array}
 \right]\right)\ge 0\ \forall x \in\bbR^m.    
\end{array}$}
\]
\label{le:sl}
\end{lemma}

From this key lemma and the fact that the ReLU $\phi:\ \bbR\to\bbR$
satisfies $\phi\in\slope(0,1)$, 
we can obtain the next result on 
the Zames-Falb multiplier for the ReLU given by \rec{eq:ReLU}.  
\begin{corollary}
Let us define
\begin{equation}
\scalebox{0.75}{$
\begin{array}{@{}l}
\arraycolsep=0.5mm
 \PiZF:=
 \left\{
  \Pi\in\bbS^{2m}:\ 
  \Pi=
 (\ast)^T
\left[
 \begin{array}{cc}
  0 & M^T\\
  M & 0 
 \end{array}
\right]
 \left[
 \begin{array}{cc}
  I_m & -I_m\\
  0 & I_m\\
 \end{array}
 \right],\ M\in \DHD^m \right\}.  
\end{array}$}
\label{eq:MZF}
\end{equation}
Then, $\Pi\in \PiZF$ is a valid multiplier that satisfies \rec{eq:IQC0}
for the ReLU $\Phi$ given by \rec{eq:ReLU}.  
\label{cor:ZF}
\end{corollary}
%

\subsection{Polytopic Bounding Multiplier}

The polytopic bounding multipliers are useful
to capture the properties of sector-bounded nonlinearities.  
To represent them in compact fashion, let us define
\begin{equation}
\scalebox{0.85}{$
\begin{array}{@{}l}
\arraycolsep=0.1mm
\Pispol[\alpha,\beta]:=
 \left\{
  \Pi \in \bbS^{2m}:\ 
  (\ast)^T \Pi 
  \left[
   \begin{array}{c}
    I \\
    \Delta 
   \end{array}
  \right]\succ 0\ \forall \Delta \in \bbD^n[\alpha,\beta]
 \right\}.  
\end{array}$}
\label{eq:Pipol0}
\end{equation}
Then the following lemma provides
the polytopic bounding multipliers for sector-bounded nonlinearities.  
\begin{lemma}\cite{Fetzler_IFAC2017}
For a given nonlinearity $\phi\in \sec[\alpha,\beta]$ with 
$\alpha\le 0\le \beta$, 
let us define $\Phi:\ \bbR^m\to \bbR^m$
by the first equation in \rec{eq:ReLU}.  
Assume $\Pi\in\Pispol[\alpha,\beta]$.  
Then we have
\[
 (\ast)^T
\Pi
\left[
 \begin{array}{c}
  x\\
  \Phi(x)
 \end{array}
 \right]\ge 0\ \forall x \in\bbR^m.    
\]
\label{le:sb}
\end{lemma}

As also stated in \cite{Fetzler_IFAC2017}, 
it is hard to check whether $\Pi\in\Pispol[\alpha,\beta]$ holds 
since $\Pispol[\alpha,\beta]$ is characterized by infinitely many constraints.  
To get around this difficulty, we employ a primitive but numerically tractable
inner approximation of $\Pispol[\alpha,\beta]$ given as follows:
\begin{equation}
\scalebox{0.9}{$
\begin{array}{@{}l}
\arraycolsep=0.3mm
\Pipol[\alpha,\beta]:=
 \left\{
  \Pi = 
  \left[
   \begin{array}{cc}
    X & Y \\
    Y^T & Z
   \end{array}
  \right]
  \in \bbS^{2m}:\right.\\ \left.
\arraycolsep=0.3mm
  (\ast)^T \Pi 
  \left[
   \begin{array}{c}
    I \\
    \Delta 
   \end{array}
  \right]\succ 0\ \forall \Delta \in \bbD_\mathrm{ver}^n[\alpha,\beta],\ 
  Z_{ii}\le 0\ (i=1,\cdots,m)
 \right\}.  \hspace*{-20mm}
\end{array}$}
\label{eq:Pipol}
\end{equation}
From this inner approximation and the fact that 
the ReLU $\phi:\ \bbR\to\bbR$
satisfies $\phi\in\sec[0,1]$, 
we can obtain the next result that provides
the polytopic bounding multiplier for the ReLU given by \rec{eq:ReLU}.  
\begin{corollary}
Let us define
\begin{equation}
\scalebox{0.9}{$
\begin{array}{@{}l}
\arraycolsep=0.5mm
\Pipol:=
 \left\{
  \Pi = 
  \left[
   \begin{array}{cc}
    X & Y \\
    Y^T & Z
   \end{array}
  \right]
  \in \bbS^{2m}:\right.\\ \left.
  (\ast)^T \Pi 
  \left[
   \begin{array}{c}
    I \\
    \Delta 
   \end{array}
  \right]\succ 0\ \forall \Delta \in \bbD_\mathrm{ver}^n[0,1],\ 
  Z_{ii}\le 0\ (i=1,\cdots,m)
 \right\}.  \hspace*{-20mm}
\end{array}$}
\label{eq:Mpol}
\end{equation}
Then, $\Pi\in \Pipol$ is a valid multiplier that satisfies \rec{eq:IQC0}
for the ReLU $\Phi$ given by \rec{eq:ReLU}.  
\label{cor:pol}
\end{corollary}

We finally note that the denomination ``polytopic bounding'' 
comes from the historical reason that 
the multipliers in \rec{eq:Pipol0} and \rec{eq:Pipol}
have been used to handle parametric uncertainties in polytopes 
in the context of robust control \cite{Iwasaki_IEEE1998,Scherer_Automatica2001}.  

\begin{remark}
Even though we restrict our attention
to the static Zames-Falb multiplier of the form \rec{eq:MZF} in \rco{cor:ZF}, 
it is true that the dynamical 
finite impulse response (FIR) Zames-Falb multipliers 
are also investigated in \cite{Fetzler_IFAC2017} in frequency domain.  
We do not pursue such a direction in this paper
mainly because the novel copositive multipliers, 
to be introduced in the next subsection, 
rely on the analysis in time-domain.  
However, we have a prospect that 
the extension similar to the FIR multipliers in 
\cite{Fetzler_IFAC2017} can also be achieved in time-domain
by means of discrete-time system lifting \cite{Bittanti_1996}.  
Such an extension, and mutual relationship with the FIR multipliers
are currently under investigation.  
Still, we have already obtained related results on the use
of the discrete-time system lifting in \cite{Ebihara_EJC2021}.  
\end{remark}

\begin{remark}
As clarified exhaustively in \cite{Fetzler_IFAC2017}, 
the polytopic bounding multiplier encompasses
some existing and frequently used multipliers.  
For instance, 
the following so-called diagonally structured multiplier has been often employed
to handle sector-bounded nonlinearities $\Phi$ in \rle{le:sb}.  
\[
 \scalebox{0.85}{$
\begin{array}{@{}l}
\bPi_\mathrm{ds}[\alpha,\beta]:=
 \left\{
  \Pi \in\bbS^{2m}:\ \Pi= 
  \left[
   \begin{array}{cc}
    -\alpha\beta D & \dfrac{\alpha+\beta}{2} D \\
    \ast & -D
   \end{array}
  \right],\ 
  D\in\bbD_{++}^m
 \right\}.  
\end{array}$}
\]
Then it is very clear that 
 $\bPi_\mathrm{ds}[\alpha,\beta]\subset\Pipol[\alpha,\beta]\subset\Pispol[\alpha,\beta]$.
 Since the effectiveness of the Zames-Falb multipliers is also widely
 recognized,
 we could say that $\Pi\in \Pipol + \PiZF$ is the most up-to-date, 
effective, and numerically tractable existing (static) multiplier to
 handle the ReLU.  
\end{remark}

\subsection{Novel Copositive Multiplier}

It has been shown recently in \cite{Raghunathan_NIPS2018}
that the input-output relationship of the ReLU given by
\rec{eq:ReLU} can be fully captured by 
three (in)equalities.  
Similar observation can also be found in \cite{Groff_CDC2019}.  
Namely, $\zeta=\phi(\xi)$ holds 
for the input $\xi\in\bbR$ and output $\zeta\in\bbR$ 
of the ReLU if and only if
\begin{equation}
\zeta(\zeta-\xi)=0,\ \zeta\ge 0,\ \zeta-\xi\ge 0.     
\label{eq:ReLU2}
\end{equation}
The first constraint is quadratic on the input and output signals
and hence compatible with IQCs.
In fact, this constraint can be regarded as the extreme case of
the sector bounded nonlinearity $\phi\in\sec[0,1]$.
From this constraint, 
we can also ensure that $\Pi\in \Pipol$ is a valid multiplier satisfying
\rec{eq:IQC0}.
On the other hand, the second and third constraints are
{\it linear} with respect to the input and output signals.  
Therefore they do not conform to the IQC framework 
if we merely rely on the standard positive semidefinite cone $\PSD$.
This is because the cone $\PSD$ has no
functionality to distinguish nonnegative vectors
in the quadratic form.  
To get around this difficulty, we employ copositive cone
$\COP$ and introduce the copositive multipliers.
This result is summarized in the next theorem.  
\begin{theorem}
Let us define
\begin{equation}
\scalebox{0.83}{$
\begin{array}{@{}l}
\arraycolsep=0.5mm
 \PisCOP:=
 \left\{\Pi\in \bbS^{2m}:\ \Pi=
 (\ast)^T Q 
 \left[
 \begin{array}{cc}
  -I_m & I_m\\
  0 & I_m\\
 \end{array}
 \right],\ Q\in \COP_{2m}
\right\}.    
\end{array}$}
\label{eq:MCOPs}
\end{equation}
Then, $\Pi\in \PisCOP$ is a valid multiplier that satisfies \rec{eq:IQC0}
for the ReLU $\Phi$ given by \rec{eq:ReLU}.  
\label{th:COP}
\end{theorem}
\begin{remark}
 As stated in \rsec{sec:cop}, it is hard to check whether
 $Q\in \COP_{2m}$ holds in \rec{eq:MCOPs} and hence
 the copositive multiplier \rec{eq:MCOPs} is intractable in general.
 To get around this difficulty, we apply inner approximation to the
 copositive cone $\COP$ and define
\begin{equation}
\scalebox{0.72}{$
\begin{array}{@{}l}
\arraycolsep=0.5mm
 \PiCOP:=
 \left\{\Pi\in \bbS^{2m}:\ \Pi=
 (\ast)^T Q 
 \left[
 \begin{array}{cc}
  -I_m & I_m\\
  0 & I_m\\
 \end{array}
 \right],\ Q\in \PSD_{2m}+\NN_{2m}
\right\}.    
\end{array}$}
\label{eq:MCOP}
\end{equation}
 Then, it is clear from \rec{eq:inc1} that $\PiCOP\subset \PisCOP$
 and hence $\Pi\in\PiCOP$ is a valid multiplier that satisfies \rec{eq:IQC0}
 for the ReLU $\Phi$ given by \rec{eq:ReLU}.
 In particular, $\PiCOP=\PisCOP$ holds if $m\le 2$.
 It should be noted that checking  $Q\in \PSD_{2m}+\NN_{2m}$ is
 numerically tractable since this is essentially a positive semidefinite
 constraint.  
\end{remark}
 \begin{remark}
In relation to the copositive multiplier \rec{eq:MCOP},
let us consider its special class given by
\[
\scalebox{0.73}{$
\begin{array}{@{}l}
\arraycolsep=0.5mm
 \PiCOPz:=
 \left\{\Pi\in \bbS^{2m}:\ \Pi=
  (\ast)^T
  \left[
   \begin{array}{cc}
    0 & 0 \\
    0 & \hatQ
   \end{array}
  \right]
 \left[
 \begin{array}{cc}
  -I_m & I_m\\
  0 & I_m\\
 \end{array}
 \right],\ \hatQ\in \PSD_{m}+\NN_{m}
\right\}.    
\end{array}$}
\]
Then, we can see from \cite{Ebihara_EJC2021} that
the condition \rec{eq:RNNIQC} with $\Pi\in \PiCOPz$
is a sufficient condition for the $l_{2+}$-induced-norm-based 
scaled small gain condition $\|D^{-1}G_0D\|_{2+}< 1$ with $D\in\bbD_{++}^m$.
Since the ReLU only returns nonnegative signals,
we intuitively deduce that $\|D^{-1}G_0D\|_{2+}< 1$ could be a
sufficient condition for the stability.
We have validated this as the main result in \cite{Ebihara_EJC2021},
providing also the numerically verifiable condition
\rec{eq:RNNIQC} with $\Pi\in \PiCOPz$.
Since $\PiCOPz\subset \PiCOP$ does hold,
we can conclude that 
the present result encompasses the main result of
\cite{Ebihara_EJC2021}
as a special case.  
\end{remark}
\begin{remark}
The treatment of nonnegative signals is the
core for the analysis of positive systems,
and to acitively use the nonnegativity in the analysis
the integral {\it linear} constraints are introduced in \cite{Briat_IJRN2013}.
However, to build an effective stability analysis method of RNNs upon  
the powerful IQC approach with existing multipliers,
we have to capture the nonnegativity of the signals in {\it quadratic} form.  
This is the reason why we introduced copositive multipliers.  
\end{remark}
%

\section{Numerical Examples}

In \rec{eq:RNN}, let us consider the case 
$\Lambda=0$, $\Wout=I_6$ and 
\[
\scalebox{0.9}{$
\begin{array}{@{}l}
\Win= 
\left[
\begin{array}{rrrrrr}
 0.29 & -0.04 &  0.02+a & -0.35 & -0.05 & -0.12\\
-0.29 & -0.24 & -0.01 &  0.12 & -0.13 &  0.18\\
-0.50 &  b &  0.23 &  0.40 & -0.28 & -0.08\\
 0.14 & -0.27 & -0.15 &  0.13 & -0.47 & -0.28\\
-0.10 & -0.10 &  0.08 &  0.14 & -0.22 &  0.50\\
-0.11 & -0.28 & -0.21 & -0.14 & -0.09 &  0.20\\
\end{array}
\right].  
\end{array}$}
\]
For $(a,b)=(0,0)$, we see $\|G_0\|_2=0.9605$.  
Here we examined the finite gain $l_2$ stability 
over the (time-invariant) parameter variation 
$a\in[-2,2]$ and $b\in [-10,10]$.  
This example is exactly the same as that of \cite{Ebihara_EJC2021}
except for the range of the parameter variation.  

We tested the following stability conditions: \\  
\ul{Test I (SSG):}\
Find $P\in \PSD_n$, $S\in\bbD_{++}^m$ such that
\rec{eq:RNNIQC} holds with $\Pi=0$.  \\
\ul{Test II ($l_{2+}$-SSG):}\ 
Find $P\in \PSD_n$, $S\in\bbD_{++}^m$, and
$\Pi\in \PiCOPz$ such that \rec{eq:RNNIQC} holds.  \\
\ul{Test III (SSG+ZF+PolB):}\ 
Find $P\in \PSD_n$, $S\in\bbD_{++}^m$ and
$\Pi\in \PiZF+\Pipol$ such that \rec{eq:RNNIQC} holds.  \\
\ul{Test IV (SSG+ZF+PolB+COP):}\ 
Find $P\in \PSD_n$, $S\in\bbD_{++}^m$ and
$\Pi\in \PiZF+\Pipol+\PiCOP$ such that \rec{eq:RNNIQC} holds.  

It is very clear that if Test I is feasible then Tests II and III are,
and if Test III is feasible then Test IV is.  
However, there is no theoretical inclusion relationship between Test II and
Test III.
Test I corresponds to the scaled small gain condition with the standard
$l_2$ induced norm,
while
Test II corresponds to the scaled small gain condition with
the $l_{2+}$ induced norm.
These have been already implemented in \cite{Ebihara_EJC2021},
but we retested them since we changed the range of the parameter variation.  

\begin{figure}[t]  
\begin{center}
 \vspace*{-2mm}
 \includegraphics[scale=0.50]{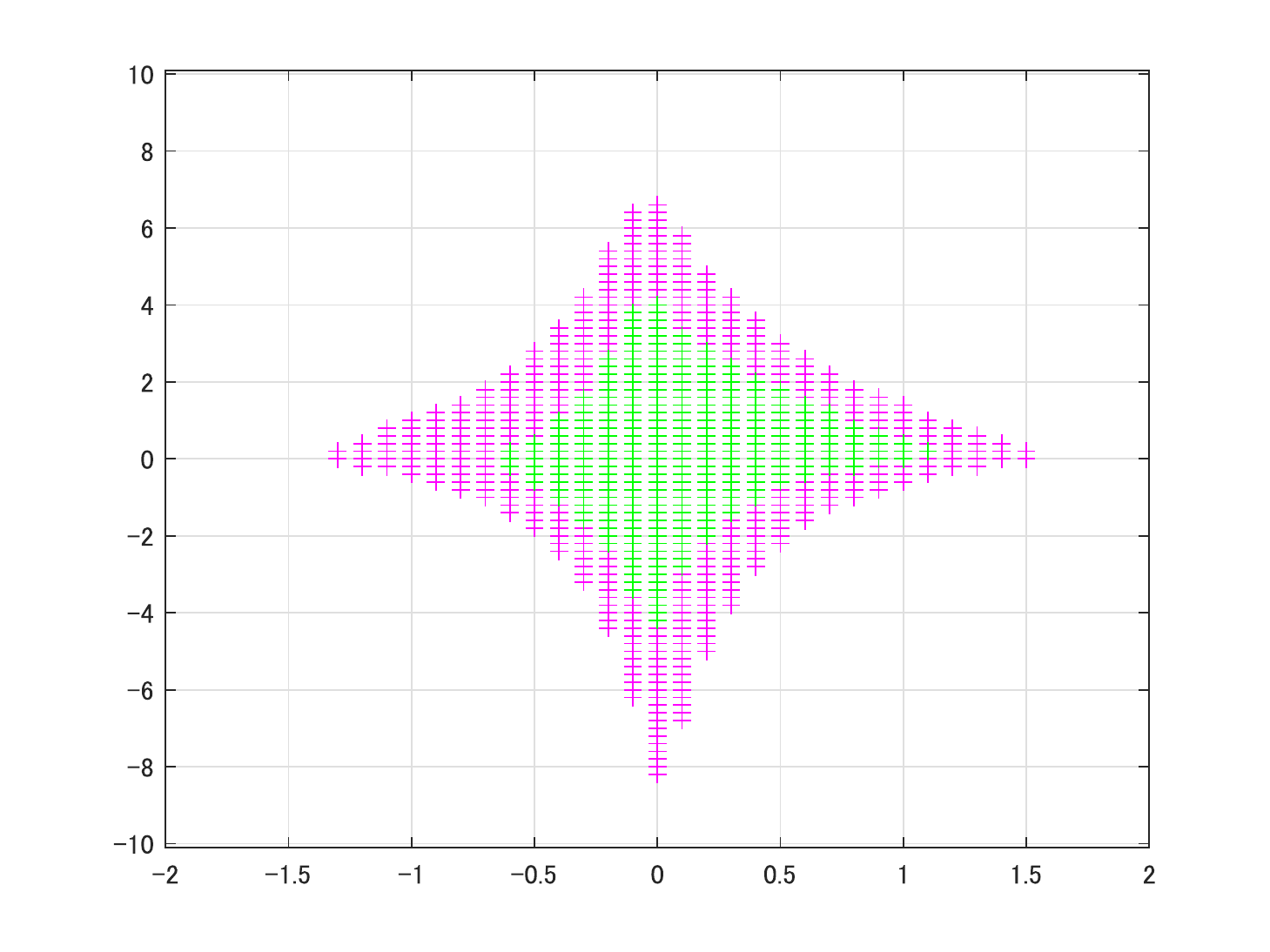}
 \vspace*{-6mm}
 \caption{Comparison: Test I  vs Test II.  }
 \label{fig:comp1}
\end{center}
 \vspace*{-4mm}
\begin{center}
 \includegraphics[scale=0.50]{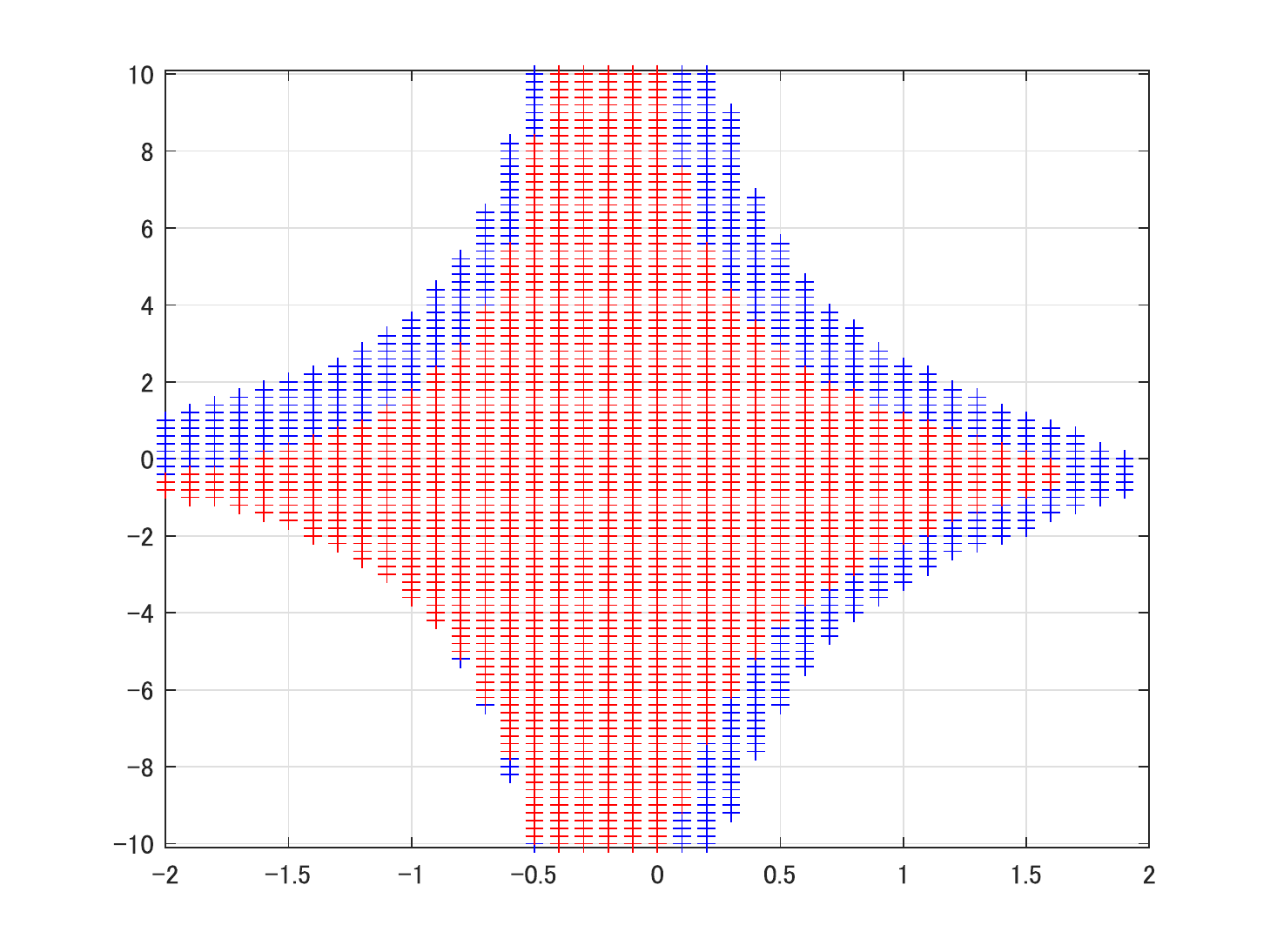}
 \vspace*{-6mm}
 \caption{Comparison: Test III  vs Test IV.  }
 \label{fig:comp2}
 \vspace*{-8mm}
\end{center}
\end{figure}

In \rfig{fig:comp1}, 
we plot $(a,b)$ for which the RNN is proved to be stable
by Tests I and II.  
Both Tests turned out to be feasible 
for $(a,b)$ in the green region, 
whereas only Test II turned out to be feasible
for $(a,b)$ in the magenta region.   
On the other hand, in \rfig{fig:comp2}, 
both Tests III and IV turned out to be feasible 
for $(a,b)$ in the red region, 
whereas only Test IV turned out to be feasible
for $(a,b)$ in the blue region.   
From both figures, we can confirm the effectiveness of the
copositive multipliers.
As for the comparison between Tests II and III,
Test III tuned out to be feasible in much larger region than that of
Test II, but there is no strict inclusion relationship between them.
In fact, for $(a,b)=(1.0,1.4)$,
Test II and III turned out to be feasible and infeasible, respectively.  

\section{Conclusion and Future Works}

In this paper, we dealt with the stability analysis of the RNN with the ReLU
by means of the IQC framework.
By actively using the nonnegativity property of the ReLU,
we newly introduced the copositive multipliers.
We showed that we can employ copositive multipliers (or their inner approximation)
together with existing multipliers such as
Zames-Falb multipliers and polytopic bounding multipliers, 
and this directly enabled us to ensure that the introduction of 
copositive multipliers leads to better (no more conservative) results.  
By numerical examples, we illustrated the effectiveness of the
copositive multipliers.  

In the present paper and \cite{Ebihara_EJC2021},
we converted a COP to an SDP
by simply replacing $\COP$ by $\PSD+\NN$.  
However, this treatment is conservative.  In this respect, 
Lasserre \cite{Lasserre_2014} and
Klerk and Pasechnik \cite{De_SIAM2002}
have already shown independently how to
construct a hierarchy of SDPs to solve COP 
in an asymptotically exact fashion, but   
the size of SDPs grows very rapidly.
This is prohibitive to deal with realistic, larger size networks.  
To get around this difficulty, 
we plan to rely on efficient first-order methods to 
solve the specific conic relaxations arising from 
polynomial optimization problems with sphere constraints
\cite{Mai_arXiv2020}.

%


\end{document}